\documentstyle{article}

\newcommand{\Frac}[2]{\displaystyle\frac{#1}{#2}}
\newcommand{\Sum}{\displaystyle\sum}
\newcommand{\N}{\mbox{\rm I\kern-1.5pt N}}
\newcommand{\F}{\mbox{${\rm I} \! {\rm F}$}}
\newcommand{\Ff}{\begin{tiny}{\rm I} \! {\rm F}\end{tiny}}
\newcommand{\PP}{\mbox{${\rm I} \! {\rm P}$}}
\newcommand{\Z}{\mbox{\sf Z \hspace{-1.1em} Z}}
\newcommand{\dem}{\par{\em Proof\/: }\\ \noindent }
\newcommand{\findemo}{$\;\;\Box$\\}

\newtheorem{defi}{Definition}[section]
\newtheorem{prop}[defi]{Proposition}
\newtheorem{thm}[defi]{Theorem}
\newtheorem{coro}[defi]{Corollary}
\newtheorem{nota}[defi]{Remark}
\newtheorem{ejplo}[defi]{Example}
\newtheorem{algor}[defi]{Algorithm}

\begin{document}

\title{Computing Weierstrass semigroups and the Feng-Rao distance 
from singular plane models}
\author{A. Campillo$\;\,^{\ast}$ \and J. I. Farr\'{a}n 
\thanks{Both authors are partially supported by DIGICYT PB97-0471}}
\date{July 25, 1999}
\maketitle

\begin{abstract}

We present an algorithm to compute 
the Weierstrass semigroup at a point $P$ 
together with functions for each value in the semigroup, 
provided $P$ is the only branch at infinity 
of a singular plane model for the curve. 
As a byproduct, the method also provides us with 
a basis for the spaces ${\cal L}(mP)$ 
and the computation of the Feng-Rao distance 
for the corresponding array of geometric Goppa codes. 
A general computation of the Feng-Rao distance is also obtained. 
Everything can be applied to the decoding problem 
by using the majority scheme of Feng and Rao. 


\end{abstract}

\section{Introduction}

One of the main problems in the theory of algebraic Geometry codes 
(AG codes in brief) 
is the explicit computation of bases for the spaces ${\cal L}(G)$, 
$G$ being a divisor over an algebraic curve. There are some general 
methods to do this computation, like Coates or Brill-Noether algorithms. 
In this paper we study an alternative procedure for the 
particular case $G=mP$, where $P$ is a rational point of the curve, 
using the structure of the Weierstrass semigroup, 
the well-known theory of Abhyankar-Moh 
and the normalization of the curve. 
This will have also the advantage that 
it can be used for effective decoding of those codes. 
Thus, our aim is to review and develop links between 
coding and singularity theories. 

Our method is algorithmic and it assumes to have a singular plane model 
for the curve with only one branch $P$ at infinity which is defined 
over the base field. This branch is nothing but the point $P$ in 
the support of $G$\/. A first step uses 
the Abhyankar-Moh theorem to give a subsemigroup of the 
Weierstrass semigroup $\Gamma$ and functions achieving its elements, 
with the aid of the so-called {\em algorithm of approximate roots}\/. 
A second step computes the rest of the semigroup and corresponding 
functions by a triangulation algorithm starting from 
an integral basis for the integral closure of the 
algebra of the affine part of the curve. 

In the case of curves over finite fields, a byproduct of our 
algorithm can be used for decoding the Algebraic Geometry codes $C(m)$ 
with divisor of type $G=mP$ by the method of Feng and Rao \cite{FR}. 
This method, based on a majority voting scheme, 
needs the knowledge of the Weierstrass semigroup at $P$ 
and functions achieving its values, and it decodes errors 
up to half the so-called Feng-Rao distance $\delta_{FR}(m)$ 
(an estimate for the minimum distance which is better than the Goppa distance).

The integer $\delta_{FR}(m)$ can be defined in terms of 
arithmetical relations among elements in the semigroup $\Gamma$\/. 
The precise value of $\delta_{FR}(m)$ is known for few classes of 
semigroups as, for instance, for (many elements of) telescopic 
semigroups (see \cite{HPNew}). Telescopic semigroups are 
complete intersection, and so a very special kind of 
symmetric semigroups (i.e. those such that $c=2g$\/, $c$ being 
the conductor and $g$ the number of gaps). 
The Abhyankar-Moh semigroup is telescopic, 
but the Weierstrass semigroup can be quite general. 

The last part of the paper is devoted to the computation of 
the Feng-Rao distance for numerical semigroups $S$\/. 
Our computational method assumes that $S$ is presented 
by means of its so-called Ap\'{e}ry systems. 
Such systems are the natural way to describe semigroups 
when one is dealing with problems involving relations 
(see \cite{CampG}). 
We show a formula (theorem \ref{formFRd}) 
to compute $\delta_{FR}(m)$ for a general $S$\/. 
If $S$ is symmetric, this formula is improved 
(theorem \ref{NmFRdSym}) for some elements with $m\geq c$\/. 
Also for symmetric semigroups, we show how the formula 
$$\delta_{FR}(m)=\min\{r\in S\;|\;r\geq m+1-2g\}\hspace{1cm}[@]$$
holds for most values $m\in S$ with $c\geq m\geq 2c-1$ and, 
moreover, a formula for the minimum element $m_{0}\in S$ 
such that the formula $[@]$ holds for $m>m_{0}\,$. 
An estimate for $m_{0}$ was given for telescopic semigroups in \cite{HPNew}. 

Finally, we show how our algorithmic method to compute 
Weierstrass semigroups $\Gamma$ gives also, as a byproduct, 
an easy way to compute an Ap\'{e}ry system for such semigroups 
and, hence, the Feng-Rao distance for their elements.

\section{Weierstrass semigroups and AG codes}

The following section outlines the connection between 
coding theory and Weierstrass semigroups, 
and it is abstracted from \cite{HohPel} and \cite{KirPel}. 
Consider a non-singular projective algebraic 
curve $\chi$ defined over a finite field $\F$ 
such that $\chi$ is irreducible over $\overline{\F}$. 
One-point geometric Goppa codes are constructed as follows. 
Take $n+1$ different $\F$\/-rational points $P_{1},\ldots,P_{n}$ 
and $P$ of the curve $\chi$\/, and take a positive integer $m$\/. 
Thus, one can consider the linear map 
$$\begin{array}{c}
ev_{D}\;:\;{\cal L}(mP)\longrightarrow\F^{n}\\
f\mapsto(f(P_{1}),\ldots,f(P_{n}))
\end{array}$$

\vspace{.1cm}

\noindent
and define the linear code $C(m)$ as the dual space of 
$Im\,(ev_{D})$\/, whose length is obviously $n$\/. 
Denote by $k(m)$ and $d(m)$ the dimension over $\F$ and 
the minimum distance of the linear code $C(m)$ respectively, 
where the integer $d(m)$ is the minimum value of non-zero entries 
of a non-zero vector of $C(m)$\/. Goppa estimates for $k(m)$ and 
$d(m)$ are derived from the {\em Riemann-Roch formula} as follows 
(see \cite{TsfVla} for more details). If $2g-2<m<n$\/, then 
$$\left\{\begin{array}{ccl}
k(m)&=&n-m+g-1\\
d(m)&\geq&m+2-2g\doteq d^{\ast}(m)
\end{array}\right.$$
where $g$ is the genus of the curve $\chi$ and $d^{\ast}(m)$ 
is the so-called Goppa designed minimum distance of $C(m)$\/. 
Apart from the excellent asymptotical behaviour of their parameters 
(see \cite{TsfVla}, for instance), the main interest of these codes 
is that they can be decoded efficiently by the majority scheme of 
Feng and Rao algorithm, which will be briefly described below. 

Fix a rational function $f_{i}\in\F(\chi)$ with only one pole at $P$ 
of order $i$ for those values of $i$ for which it is possible, 
i.e. for the non-negative integers in the Weierstrass semigroup 
$\Gamma=\Gamma_{P}$ of $\chi$ at $P$\/. 
For a received word $\mbox{\bf y}=\mbox{\bf c}+\mbox{\bf e}$, 
where $\mbox{\bf c}\in C(m)$, one can consider the unidimensional 
and bidimensional syndromes given respectively by 
$$s_{i}(\mbox{\bf y})\doteq\sum_{k=1}^{n}e_{k}\,f_{i}(P_{k})
\;\;\;\;\;\;{\rm and}\;\;\;\;\;\;
s_{i,j}(\mbox{\bf y})\doteq\sum_{k=1}^{n}e_{k}\,f_{i}(P_{k})\,f_{j}(P_{k})$$
 
\noindent
Notice that the set $\{f_{i}\;|\;i\leq m,\;\; i\in\Gamma\}$ 
is actually a basis for ${\cal L}(mP)$ and hence one has 
$$C(m)=\{\mbox{\bf y}\in\F^{n}\;|\;s_{i}(\mbox{\bf y})=0\;\;
{\rm for}\;\;i\leq m\}$$
Therefore one can calculate $s_{i}(\mbox{\bf y})$ 
from the received word ${\bf y}$ as 
$$s_{i}(\mbox{\bf y})=\sum_{k=1}^{n}y_{k}\,f_{i}(P_{k})$$
and thus $s_{i}(\mbox{\bf y})$ is called {\em known}\/ for $i\leq m$\/. 

In fact, it is a known fact that if one had enough 
syndromes $s_{i,j}(\mbox{\bf y})$ for $i+j>m$ 
one could know the emitted word $\mbox{\bf c}$, 
and all the involved syndromes can be computed 
by a majority voting (see \cite{FR}). 

This decoding algorithm corrects up to half the so-called 
Feng-Rao distance $\delta_{FR}(m)$ of $C(m)$ (see definition \ref{FRdist} below).  
More precisely, the procedure corrects up to 
$\left\lfloor\Frac{\delta_{FR}(m')-1}{2}\right\rfloor$ errors, 
where $m'=\min\{r\in\Gamma\;|\;r>m\}$\/. 
In particular, one has $d(m)\geq\delta_{FR}(m')$ 
(see \cite{KirPel} for a direct proof of this inequality). 

Notice that $\delta_{FR}(m)$ only depends on the semigroup $\Gamma$\/. 
On the other hand, it is known that 
$\delta_{FR}(m)\geq m+1-2g=d^{\ast}(m-1)$ for all $m\in\Gamma$\/, 
the right hand side term being the Goppa estimate for the minimum distance. 
Thus, the method gives an improvement of the number of errors that one can 
correct for one-point AG codes. 

In practice, the main problem is computing $\Gamma$ and the 
functions $f_{i}$ achieving the values of the semigroup $\Gamma$ 
in order to carry out this decoding algorithm. On the other hand, 
a second basic problem is calculating the value of $\delta_{FR}(m)$ if 
one wants to know the number of errors that one can correct. 
The aim of this paper is to give an approach to the above two 
problems with certain conditions.

\section{Computing Weierstrass semigroups}

In this section we show how the Weierstrass semigroups 
of curves having a singular plane model with only one 
branch at infinity can be computed at the same time 
as functions achieving their values. Everything 
is based on the theory of Abhyankar and Moh about such 
plane models, and so this section can be regarded as an 
application of singularity theory to coding theory.

\subsection{Semigroups of plane curves with only one branch at infinity}

In the sequel, we will assume that $\F$ is any perfect field, 
since the hypothesis on the finiteness of $\F$ will not be used. 
Nevertheless, we are always thinking of finite fields, 
because of the applications in coding theory. 
Let $\chi$ be a non-singular projective algebraic curve defined 
over $\F$ such that $\chi$ is irreducible over $\overline{\F}$. 
Let $\Upsilon$ be a plane model for $\chi$ with only one branch at infinity, 
i.e. such that there exist a birational morphism 
$${\rm n}\;:\;\chi\rightarrow \Upsilon\subseteq \PP^{2}$$
and a line $L\subseteq\PP^{2}$ defined over $\F$ such that 
$L\cap \Upsilon$ consists of only one geometric point $P$ 
and $\chi$ has only one branch at $P$\/. 
Notice that, a fortiori, both $P$ and the branch at $P$ are 
defined over the underlying field $\F$, since $\Upsilon$ does. 
Thus there is only one point of $\chi$ lying over $P$, which will be 
denoted by $\overline{P}$\/. 

Set ${\cal C}=\chi\setminus\{\overline{P}\}$ 
and ${\cal Z}=\Upsilon\setminus\{P\}$\/. 
One has the two following additive subsemigroups of $\N$: 
$$
\begin{array}{l}
\Gamma_{P}\doteq\{-\upsilon_{\overline{P}}(f)\;|\;
f\in{\cal O}_{\chi}({\cal C})\}\\
S_{P}\doteq\{-\upsilon_{\overline{P}}(f)\;|\;
f\in{\cal O}_{\Upsilon}({\cal Z})\}
\end{array}
$$

The first one is nothing but the Weierstrass semigroup 
of $\chi$ at $\overline{P}$; 
this semigroup contains the second one but they are different 
unless the curve $\Upsilon$ is non-singular in the affine part. 
Notice that both $\N\setminus\Gamma_{P}$ and $\N\setminus S_{P}$ are finite. 
In fact $\N\setminus\Gamma_{P}$ has $g$ elements, 
$g$ being the genus of $\chi$, which are the Weierstrass gaps. 
In order to compute the cardinality of $\N\setminus S_{P}$, we need 
the following effective version of a standard fact in singularity theory.

\begin{prop}\label{TriangLem}

Let $A={\cal O}_{\chi}({\cal C})\;$ and $\;B={\cal O}_{\Upsilon}({\cal Z})$ 
be the respective affine coordinate $\F$\/-algebras of the curves ${\cal C}$ 
and ${\cal Z}$\/; then one has 
$$\sharp(\Gamma_{P}\setminus S_{P})=dim_{\Ff}(A/B)
=\sum_{Q\in{\cal Z}}\delta_{Q}({\cal Z})$$
where $Q$ ranges over all the closed points of ${\cal Z}$ 
and $\delta_{Q}({\cal Z})=
dim_{\Ff}(\overline{\cal O}_{\Upsilon,Q}/{\cal O}_{\Upsilon,Q})$, 
$\overline{\cal O}_{\Upsilon,Q}$ being the normalization of the ring 
${\cal O}_{\Upsilon,Q}\,$. 

\end{prop}

\dem

The second equality follows from the fact that 
$$A/B=\bigoplus_{Q\in Sing\,{\cal Z}}
\overline{\cal O}_{\Upsilon,Q}/{\cal O}_{\Upsilon,Q}$$ 
where $Sing\,{\cal Z}$ is the set of (closed) singular points of ${\cal Z}$\/. 
In order to prove the first equality, we take 
an $\F$\/-basis $\{h_{1},\ldots,h_{s}\}$ of $A/B$\/, 
which can be computed with the aid of the {\em integral basis algorithm}\/ 
\footnote{After \cite{Hoe}, the {\em integral basis algorithm}\/ 
is efficient whenever the ramification of $B|\F$ is not wild. 
For our purpose, this procedure will be applied in the hypothesis of proposition 
\ref{AMth}, what actually implies that the ramification is tame. } 
(see \cite{Hoe}, \cite{Trag} and also \cite{Thes}). 

Set $B^{i}\doteq B+\F h_{1}+\ldots+\F h_{i}$, for $0\leq i\leq s$\/. 
We proceed by induction, so let $0\leq i<s$ 
and assume that we have found functions $g_{1},\ldots,g_{i}$ 
which are linearly independent over $\F$ such that 
$$B+\F g_{1}+\ldots+\F g_{i}=B^{i}$$
and 
$$-\upsilon_{\overline{P}}(g_{j})\notin\Gamma_{P}^{j-1}\;\;\;\;\forall j\leq i$$
where $\Gamma_{P}^{j}\doteq S_{P}\cup\{-\upsilon_{\overline{P}}(g_{1}),\ldots,
-\upsilon_{\overline{P}}(g_{j})\}\subseteq\Gamma_{P}\,$. 

Now look at $h_{i+1}\,$. 
If $-\upsilon_{\overline{P}}(h_{i+1})\notin\Gamma_{P}^{i}$, 
then set $g_{i+1}=h_{i+1}$ and go on. 

Otherwise, there exists $f\in B^{i}$ with 
$\upsilon_{\overline{P}}(h_{i+1})=\upsilon_{\overline{P}}(f)$ and such that 
$$-\upsilon_{\overline{P}}(h_{i+1}-f)<-\upsilon_{\overline{P}}(h_{i+1})$$
Then we repeat the process with $h_{i+1}-f$ replacing to $h_{i+1}\,$. 
Since one has $h_{i+1}\notin B^{i}$, it follows that 
in a finite number of steps we will be able to replace $h_{i+1}$ by 
$$g_{i+1}\equiv h_{i+1}\;\;(mod\;\;B^{i})$$ 
such that 
$$-\upsilon_{\overline{P}}(g_{i+1})\notin\Gamma_{P}^{i}$$

\newpage

At the end of the procedure $s$ different elements in 
$\Gamma_{P}\setminus S_{P}$ will be added, and hence 
$\sharp(\Gamma_{P}\setminus S_{P})\geq dim_{\Ff}(A/B)$\/. 
On the other hand, since $A=B^{l}=B+\F g_{1}+\ldots+\F g_{s}$ 
one has that any $h\in A$ can be written in a (unique) way as 
$h=g+\lambda_{1}g_{1}+\ldots+\lambda_{s}g_{s}$ 
with $g\in B$ and $\lambda_{i}\in\F$, the values $\upsilon_{\overline{P}}(g_{i})$ 
and $\upsilon_{\overline{P}}(g)$ being pairwise different. 
Thus, one has either $-\upsilon_{\overline{P}}(h)\in S_{P}\;$ or 
$\;-\upsilon_{\overline{P}}(h)=-\upsilon_{\overline{P}}(g_{i})$ for a unique $i$\/. 
This proves the equality. 

\vspace{.1cm}
\findemo

Now our aim is as follows: we intend to describe $\Gamma_{P}$ by 
first computing $S_{P}$ and the associated functions, 
and secondly by completing the semigroup up to $\Gamma_{P}$ with 
the corresponding functions. According to the above proposition, 
this last thing can be done by means of the following procedure.

\begin{algor}[Triangulation algorithm]\label{TriangAlg}

$\;$

\vspace{.2cm}

\noindent
Input: $S_{P}\,$, ${\cal O}_{\Upsilon}({\cal Z})$ 
and $\{h_{1},\ldots,h_{s}\}$\/. 

\begin{itemize}

\item Initialize $S=S_{P}$ and $B={\cal O}_{\Upsilon}({\cal Z})$ 

\item For $i=1,\ldots,s$ do 

\begin{itemize}

\item Set $g_{i}=h_{i}$ 

\item While $-\upsilon_{\overline{P}}(g_{i})\in S$ do 

\begin{itemize}

\item Find $f\in B$ such that 
$-\upsilon_{\overline{P}}(g_{i}-f)<-\upsilon_{\overline{P}}(g_{i})$ 

\item Set $g_{i}=g_{i}-f$ 

\end{itemize}

\item Set $B=B+\F g_{i}$ and 
$S=S\cup\{-\upsilon_{\overline{P}}(g_{i})\}$ 

\end{itemize}

\item Next $i$ 

\end{itemize}

\noindent
Output: $S$ and $B$\/. 

\end{algor}

\begin{nota}\label{fastTA}

One can substitute the set $S=S\cup\{-\upsilon_{\overline{P}}(g_{i})\}$ 
by the semigroup $S=S-\upsilon_{\overline{P}}(g_{i})\cdot\N$ in the above 
algorithm, i.e. in each step one can cover more than one gap and 
also add the corresponding functions (products of existing data functions) 
to $B$\/. It may yield in general a faster algorithm, since one could stop 
just when $l$ new values are added. This idea can be done effective 
after the results of the next section by using Ap\'{e}ry systems 
(see remark \ref{sGsAp}). 

\end{nota}

\subsection{Approximate roots}

The remaining part of the algorithm, that is, the description of the 
semigroup $S_{P}$ and the construction of the corresponding functions, 
is known from the Abhyankar-Moh theorem and the so-called 
{\em algorithm of approximate roots} which, 
for the sake of completeness, will be restated below. 
Both were introduced and developed by Abhyankar and Moh in \cite{AbhMoh}. 
In particular, it will be shown that this way 
to compute the Weierstrass semigroup and the functions is effective 
and, moreover, the arithmetic properties of the involved semigroups 
will allow us to compute the Feng-Rao distance with the aid of some 
extra techniques, what will be done in section {\bf 4.3}. 
First we need the notion of approximate root.

Thus, let $S$ be a ring, $G\in S[{\rm Y}]$ a monic polynomial of 
degree $e$ and $F\in S[{\rm Y}]$ a monic polynomial of degree 
$m$ with $e|m$\/. If we write $m=ed$\/, then $G$ is called an 
approximate $d$\/-th root of $F$ if $deg\,(F-G^{d})<m-e=e\,(d-1)$\/. 
The key result is that, provided $d$ is a unit in the ring $S$\/, 
then there exists a unique approximate $d$\/-th root of $F$\/, 
which will be denoted $app(d,F)$\/. You can see a 
constructive proof in \cite{AbhSem}, where it is shown in particular 
that the computation of approximate roots is very efficient. 
From now on, we will work with the coefficient ring $S=\F[{\rm X}]$\/.

Let the affine plane model $\Upsilon$ having only one point at infinity 
be given by the equation 
$$F=F({\rm X},{\rm Y})={\rm Y}^{m}+a_{1}({\rm X})\,{\rm Y}^{m-1}+
\ldots+a_{m}({\rm X})$$
where $m$ is actually the total degree of the polynomial $F$\/, 
and set $n\doteq deg_{\rm X}F$\/. 

In the sequel, we will assume moreover that 
\begin{quote}
(H) $\;\;char\,\F$ does not divide either $deg\,\Upsilon$ or $e_{P}(\Upsilon)$
\end{quote}

\noindent
With the above notations, 
one has $m=deg\,\Upsilon$ and $n=deg\,\Upsilon-e_{P}(\Upsilon)$\/, and thus 
(H) is equivalent to say that $p=char\,\F$ does not divide either $m$ or $n$\/. 

In this case, we can actually assume that $p$ does not divide $m=deg\,\Upsilon$\/. 
In fact, if $m$ is a multiple of $p$ but $n$ is not, we choose $k$ not divisible 
by $p$ such that $nk>m$\/, and by doing a change of variables in the form 
${\rm X}'={\rm X}+{\rm Y}^{k}$, ${\rm Y}'={\rm Y}$ we get a new affine curve 
which is isomorphic to the original one\footnote{In particular, $S_{P}$ does not change. } 
but whose degree is not divisible by $p$\/. 
This will be assumed in the sequel for simplicity, 
and thus we will be able to compute any approximate root of $F$\/.

By the sake of economy one wants to avoid the computation of parametric 
equations for the singularity at infinity, and thus computing the 
semigroup $S_{P}$ and corresponding functions directly from the equation 
of the curve. The parametrization would allow us to compute easily 
$\upsilon_{\overline{P}}$\/-values, but this can be avoided by using 
resultants and approximate roots, as it is shown below.

In order to state the {\em algorithm of approximate roots} 
from \cite{AbhSem}, we first agree to set 
$$deg_{\rm X}Res_{\rm Y}(G,H)=-\infty\;\;{\rm if}\;\;
Res_{\rm Y}(G,H)=0$$
for any couple of polynomials $G,H\in\F[{\rm X},{\rm Y}]$\/, and
$$g\,c\,d\;(\delta_{0},\delta_{1},\ldots,\delta_{i})=
g\,c\,d\;(\delta_{0},\delta_{1},\ldots,\delta_{j})$$
if $\delta_{0},\delta_{1},\ldots,\delta_{j}$ are integers, $j<i$ and 
$\delta_{j+1}=\delta_{j+2}=\ldots=\delta_{i}=-\infty$\/. 
Thus, the algorithm works as follows with an input $F$ as above 
(the case when ${\rm Y}$ divides $F$ being trivial, we assume the contrary).

\begin{algor}[Approximate roots]\label{AlgAR}

$\;$

\begin{itemize}

\item Set $d_{0}=0$, $F_{0}={\rm X}$, $\delta_{0}=d_{1}=m$\/, 
$F_{1}={\rm Y}$ and $\delta_{1}=deg_{\rm X}Res_{\rm Y}(F,F_{1})$

\item For $i$ from $2$ do 

\begin{itemize}

\item $d_{i}=g\,c\,d\;(d_{i-1},\delta_{i-1})$ 

\item If $d_{i}=d_{i-1}$ then $h=i-2$ and STOP else 

\begin{itemize}

\item $F_{i}=app(d_{i},F)$ 

\item $\delta_{i}=deg_{\rm X}Res_{\rm Y}(F,F_{i})$ 

\end{itemize}

\end{itemize}

\item Next $i$ 

\end{itemize}

\noindent
Output: $h$\/, $(\delta_{0},\ldots,\delta_{h})$ and $(F_{0},\ldots,F_{h})$ 

\end{algor}

\vspace{.3cm}

Since the sequence $\{d_{i}\}_{i\geq 1}$ 
obtained in the above algorithm 
is a decreasing one of positive integers, 
there exists a unique positive integer $h$ 
such that $d_{1}>\ldots>d_{h+1}=d_{h+2}\,$, 
and hence the algorithm stops. The first application of 
algorithm \ref{AlgAR} is the following result, proved by Abhyankar in \cite{AbhIrr}, 
which provides a criterion for a curve with an only (rational) point 
at infinity to have only one (rational) branch at this point 
(and to be absolutely irreducible, as a consequence).

\begin{prop}[{\bf Criterion for having only one branch at infinity}]\label{criter}

$\;$

\noindent
Let $F$ be the equation of a plane model with an only point at 
infinity as above, and assume that $char\,\F$ does not divide 
$m=deg\,F$\/. Let $h$\/, $d_{i}$ and $\delta_{i}$ the integers 
which are computed by the algorithm of approximate roots. 
Then the curve has an only (rational) branch at infinity if and only if 
$d_{h+1}=1$, $\delta_{1}d_{1}>\delta_{2}d_{2}>\ldots>\delta_{h}d_{h}$ 
and $n_{i}\delta_{i}$ is in the semigroup generated by 
$\delta_{0},\delta_{1},\ldots,\delta_{i-1}$ for $1\leq i\leq h$, 
where $n_{i}\doteq d_{i}/d_{i+1}$ also for $1\leq i\leq h$\/. 

\end{prop}

The second application of the algorithm of approximate roots is 
just the computation of $S_{P}$ and the associated functions by 
means of the Abhyankar-Moh theorem, which provides us with a set 
of generators for $S_{P}$ with nice arithmetic properties. 
The proof is refered to \cite{AbhMoh} or \cite{Pink}.

\begin{prop}[Abhyankar-Moh theorem]\label{AMth}

$\;$

\noindent
Let $\Upsilon$ be a plane model with an only branch at infinity, 
and assume that $char\,\F$ does not divide $deg\,\Upsilon$\/. 
Then there exist an integer $h$ and a sequence of integers 
$\delta_{0},\ldots,\delta_{h}\in S_{P}$ generating $S_{P}$ 
such that $\delta_{0}=deg\,\Upsilon$ and 

\begin{description}

\item[(I)] $d_{h+1}=1$ and $n_{i}>1$ for $2\leq i\leq h$\/, where 
$d_{i}\doteq g\,c\,d\;(\delta_{0},\ldots,\delta_{i-1})$ for 
$1\leq i\leq h+1$ and $n_{i}\doteq d_{i}/d_{i+1}$ for 
$1\leq i\leq h$\/. 

\item[(II)] $n_{i}\delta_{i}$ is in the semigroup generated by 
$\delta_{0},\ldots,\delta_{i-1}$ for $1\leq i\leq h$ 
\footnote{Notice that $n_{1}\delta_{1}$ is always a multiple of 
$\delta_{0}\,$. Thus the three properties of the Abhyankar-Moh 
theorem are trivially satisfied if $h\leq 1$ and $d_{h+1}=1$. }. 

\item[(III)] $n_{i}\delta_{i}>\delta_{i+1}$ for $1\leq i\leq h-1$\/. 

\end{description}

\end{prop}

\vspace{.3cm}

\begin{nota}\label{counterEx}

Notice that the restriction on the characteristic of $\F$ is necessary, 
since the plane curve over $\F_{2}$ given by the equation
$${\rm Y}^{8}+{\rm Y}={\rm X}^{2}({\rm X}^{8}+{\rm X})$$
has no affine singularity and an only point at infinity 
whose Weierstrass semigroup is generated by the elements $\{8,10,12,13\}$\/. 
This example gives a negative answer to a question 
proposed by Pinkham in \cite{Pink}, since the sequences 
$\{8,12,10,13\}$ or $\{12,8,10,13\}$ actually satisfy 
{\bf (I)}, {\bf (II)} and {\bf (III)}, but $\delta_{0}=8$ or $\delta_{0}=12$, 
respectively, is not the degree of the curve. 
Thus, a general description for $\Gamma_{P}=S_{P}$ in this way 
is still an open problem. 

Curiously, this example comes from coding theory (see \cite{HanStich}), 
and thus the reason why the Abhyankar-Moh theorem fails is that both 
$m$ and $e_{P}(\Upsilon)$ are multiple of $2$. 

\end{nota}

\begin{nota}\label{sGs}

Semigroups as in proposition \ref{AMth}, 
as well as semigroups of values studied in positive characteristic by Angerm\"{u}ller 
\footnote{In such semigroups, property {\bf (III)} in proposition \ref{AMth} 
is substituted by $n_{i}\delta_{i}<\delta_{i+1}\,$. }
in \cite{Ang} or Campillo in \cite{Camp}, are a particular case of telescopic 
semigroups, where only the properties {\bf (II)} and $d_{h+1}=1$ are required 
(see \cite{KirPel}). 
Being telescopic is equivalent to be free in the sense of \cite{Ang}, what means 
that every $m\in S_{P}$ can be written in a unique (and effective) way in the form 
$$m=\Sum_{i=0}^{h}\lambda_{i}\delta_{i}\;\;\;\;\;\;\;\;[\star]$$
with $\lambda_{0}\geq 0$ and $0\leq\lambda_{i}<n_{i}$ 
for $1\leq i\leq h$ (see \cite{KirPel}). 
From property {\bf (II)} one has that these semigroups are complete intersection 
(i.e. those such that the affine curve defined by them is a complete intersection one)
and, in particular, they are symmetric (see \cite{HerzKunz}, \cite{KirPel} 
and section {\bf 4.1} below for further details). 

\end{nota}

The functions which are needed to achieve the pole orders given by 
the Abhyankar-Moh theorem can actually be assumed to be successive  
approximate roots of the equation $F$\/, and they can be computed 
by the algorithm of approximate roots. More precisely, one has 
$$\delta_{i}=I_{P}(F,F_{i})=deg_{\rm X}Res_{\rm Y}(F,F_{i})
\;\;\;\;{\rm for}\;\;1\leq i\leq h$$

\newpage

\noindent
where $F_{i}=app(d_{i},F)$\/, $I_{P}(F,F_{i})$ denotes the intersection 
number at $P$ of the projective completions of the curves given by $F$ and $F_{i}$\/, 
and where $F_{i}$ and $\delta_{i}$ are given by the algorithm of approximate roots 
(see \cite{AbhExp}). In particular, one has 
$-\upsilon_{\overline{P}}(F_{i})=\delta_{i}\,$. 
Thus, algorithm \ref{AlgAR} computes at the same time 
the generators provided by the Abhyankar-Moh theorem 
and functions achieving such values in an effective way.

Even more, if the algorithm succeeds, i.e. if one arrives to the end 
with the properties required by proposition \ref{criter} 
\footnote{Notice that the properties given by the Abhyankar-Moh theorem 
are just the same which are required by the criterion for one branch at infinity. }, 
we are sure that there is one branch at $P$\/, and hence the curve is absolutely irreducible. 
In case of fail, i.e. if such conditions are not satisfied in any of the steps of the algorithm, 
one can conclude that the plane curve has more than one branch at $P$\/.

\begin{nota}\label{complexity}

In case of having a priori a parametrization for the singularity of the 
unique branch at $P$\/, we could use it to compute contact orders instead 
of using resultants because of the formula $\delta_{i}=I_{P}(F,F_{i})$\/. 
This would give an alternative for the algorithm. 

\end{nota}

\begin{ejplo}\label{exAAR}

Consider the affine plane curve ${\rm Y}^{8}+{\rm Y}^{2}+{\rm X}^{3}=0$ 
defined over $\F_{2}$, with only one point at infinity $P=(1:0:0)$. 
The degree of the curve is multiple of the characteristic, so with 
the change ${\rm X}={\rm X}+{\rm Y}^{3}$, ${\rm Y}={\rm Y}$ one gets 
the plane model $F({\rm X},{\rm Y})={\rm Y}^{9}+{\rm Y}^{8}+
{\rm X}{\rm Y}^{6}+{\rm X}^{2}{\rm Y}^{3}+{\rm Y}^{2}+{\rm X}^{3}$\/, 
and one can apply the algorithm of approximate roots to $F$: 
$$F_{0}={\rm X}\, ,\,\delta_{0}=d_{1}=9\, ,\,F_{1}={\rm Y}$$
$$\delta_{1}=deg_{\rm X}Res_{\rm Y}(F,{\rm Y})=3\, ,\,
d_{2}=g\,c\,d\,(9,3)=3$$
$$F_{2}=app(3,F)={\rm Y}^{3}+{\rm Y}^{2}+{\rm Y}+{\rm X}+1$$
$$\delta_{2}=deg_{\rm X}Res_{\rm Y}(F,F_{2})=8\, ,\,
d_{3}=g\,c\,d\,(9,3,8)=1$$
Thus $h=2$ and $S_{P}=\langle 9,3,8 \rangle$. 
As a consequence, there is only one branch at infinity 
since properties {\bf (I)}, {\bf (II)} and {\bf (III)} from proposition 
\ref{AMth} are satisfied. 

On the other hand, with the notations as in proposition \ref{TriangLem}, 
take a $\F_{2}$-basis for $A/B$: 
$$h_{1}=\frac{{\rm Y}(1+{\rm Y}^{6})}{{\rm X}+{\rm Y}^{3}}\;\;\;\; 
h_{2}=\frac{{\rm Y}(1+{\rm Y}^{6})}
{({\rm X}+{\rm Y}^{3})({\rm Y}^{2}+{\rm Y}+1)}$$ 
$$h_{3}=\frac{{\rm X}^{2}+{\rm Y}^{6}}{{\rm Y}^{2}+{\rm Y}+1}\;\;\;\; 
h_{4}=\frac{{\rm Y}^{2}(1+{\rm Y}^{3})({\rm Y}^{2}+{\rm Y}+1)}
{{\rm X}+{\rm Y}^{3}}$$
The values of this functions are $-\upsilon_{P}(h_{1})=13\notin S_{P}$, 
$-\upsilon_{P}(h_{2})=7\notin\Gamma_{P}^{1}$, 
$-\upsilon_{P}(h_{3})=10\notin\Gamma_{P}^{2}$ and 
$-\upsilon_{P}(h_{4})=13\in\Gamma_{P}^{3}$. Then change $h_{4}$ by 
$$g_{4}=h_{4}+h_{1}=\frac{{\rm Y}(1+{\rm Y}^{3})({\rm Y}^{2}+{\rm Y}+1)}
{{\rm X}+{\rm Y}^{3}}$$
and now $-\upsilon_{P}(g_{4})=10\in\Gamma_{P}^{3}$. Thus, one still has 
to take the function $$g_{4}=h_{4}+h_{1}+h_{3}=
\frac{{\rm Y}(1+{\rm Y}^{3})({\rm Y}^{4}+{\rm Y}^{2}+1)+
({\rm X}+{\rm Y}^{3})^{3}}
{({\rm X}+{\rm Y}^{3})({\rm Y}^{2}+{\rm Y}+1)}$$
and now $-\upsilon_{P}(g_{4})=4\notin\Gamma_{P}^{3}$. 
Hence, the Weierstrass semigroup at $P$ is 
$$\Gamma_{P}=\{0,3,{\bf 4},6,{\bf 7},8,9,{\bf 10},11,12,{\bf 13},14,\ldots\}$$

\end{ejplo}

\subsection{Application to codes}

In particular, the above results allow us to compute 
a basis of the vector space ${\cal L}(mP)$\/, 
for every $m\in\Gamma_{P}\,$, by collecting 
just one function with an only pole at $P$ of order $r$\/, 
for each $0\leq r\leq m$ with $r\in\Gamma_{P}\,$. 
This is essential for the construction and decoding of one-point 
AG codes, and it can be easily done from the above results as follows.

If $r\in S_{P}$, one gets in an effective way the writing $[\star]$ 
from the remark \ref{sGs} in the form $r=
\lambda_{0}\delta_{0}+\lambda_{1}\delta_{1}+\ldots+\lambda_{h}\delta_{h}\,$. 
Thus, if we have functions $f_{i}$ with $\upsilon_{P}(f_{i})=-\delta_{i}$ 
(for instance, those which are obtained from the algorithm of approximate 
roots), then $f_{r}=f_{0}^{\lambda_{0}}\cdot f_{1}^{\lambda_{1}}\cdot\ldots
\cdot f_{h}^{\lambda_{h}}$ has an only pole at $P$ of order $r$\/. 
Otherwise, if $r\in\Gamma_{P}\setminus S_{P}$ 
the function $f_{r}$ is constructed by algorithm \ref{TriangAlg}.

\begin{nota}\label{Reguera}

If we fix previously a semigroup with properties as in 
proposition \ref{AMth}, we can try to find a plane curve 
$F$ with an only branch at infinity achieving this semigroup, 
even if $p$ divides $g\,c\,d\;(\delta_{0},\delta_{1})$\/. 

This last procedure is just the inverse of that we have explained in 
this section, i.e. let the numbers $\delta_{0},\delta_{1},\ldots,
\delta_{h}$ with the properties as in proposition \ref{AMth} be given, 
then one can construct in a recurrent way a sequence of polynomials 
which are called \lq\lq approximants" of a polynomial $F$ such that 
$F=0$ is a curve with a unique branch $P$ at infinity 
and semigroup $S_{P}$ generated by 
$\delta_{0},\delta_{1},\ldots,\delta_{h}$ 
(this is shown by Reguera in \cite{Reg}, where one can check more details). 
This has the advantage that one can built directly the examples 
together with functions (the approximants) having the generators of 
$S_{P}$ as pole orders. 

\end{nota}

\begin{nota}\label{improvedAG}

By using the affine algebra $B$ instead of the normal affine algebra $A$ 
(with the notations as in section {\bf 3.1}), one can construct 
\lq\lq improved Algebraic Geometry codes" as in \cite{HPNew}. 
More precisely, such codes would be the dual of the functional codes 
given by ${\rm Im}\,(ev_{D}({\cal L}(mP)\cap B))$\/. 
They can also be decoded by the Feng and Rao method in the same way. 
Now one only needs the semigroup $S_{P}$ instead of $\Gamma_{P}\,$, 
so one avoids the triangulation algorithm \ref{TriangAlg}. 
Moreover, by using remark \ref{Reguera} one can directly construct 
curves with only one branch at infinity and functions achieving 
the values of a semigroup $S_{P}$ fixed a priori. 

\end{nota}

Finally, the results of this section can be summarized as follows.

\begin{thm}\label{sum3}

Let $\Upsilon$ be an absolutely irreducible projective plane curve of degree $m$ 
with only one branch $P$ at infinity which is rational over the base field $\F$. 
Assume the characteristic of $\F$ is either $0$ or it does not divide 
simultaneously $m$ and $m-e_{P}(\Upsilon)$. 
Then, by combining algorithms \ref{TriangAlg} and \ref{AlgAR} 
one can compute the Weierstrass semigroup $\Gamma_{P}$ and 
functions achieving the pole orders in $\Gamma_{P}\,$. 
As a byproduct, one gets a basis of 
the vector space ${\cal L}(rP)$, for every $r\in\Gamma_{P}\,$. 
Moreover, for each semigroup and generators with the properties as in 
the Abhyankar-Moh theorem, one can generate curves with $S_{P}$ equal 
to that semigroup, without restrictions on the characteristic. 

\end{thm}

\findemo

In the next section, we will show how to calculate the Feng-Rao distance 
for Weierstrass semigroups which have been computed with the above method. 
In fact, we will focus on a general situation in arithmetic semigroups, 
just taking into account the arithmetic properties of the semigroups in 
the theory of Abhyankar-Moh and the modifications given by proposition \ref{TriangLem}.

\section{Computing the Feng-Rao distance}

In this section, we will compute the Feng-Rao distance as a 
function defined on any arbitrary numerical semigroup, i.e. 
a subsemigroup of $\N$. Some formulae will stand for the general case, 
and then we will focus on some concrete types of semigroups, 
c.g. symmetric, Abhyankar-Moh (or, more generally, telescopic) 
and Weierstrass semigroups obtained as semigroups at infinity. 
Notice that such types of semigroups are also interesting in singularity theory.

\subsection{Ap\'{e}ry systems and Feng-Rao distance}

In the sequel, we consider numerical semigroups, 
i.e. subsemigroups $S$ of $\N$ such that 
$\sharp(\N\setminus S)<\infty$ and $0\in S$\/. 
The number $g\doteq\sharp(\N\setminus S)$ is called the {\em genus} of the 
semigroup $S$\/. Since the genus is finite, there exists a (unique) element 
$c\in S$ such that $c-1\notin S$ and $c+l\in S$ for all $l\in\N$. The number 
$c$ is called the {\em conductor} of $S$\/, and one has $c\leq 2g$\/. 
Thus, the \lq\lq last gap" of $S$ is $l_{g}\doteq c-1\leq 2g-1$, 
where $k$ is called a gap of $S$ if $k\in\N\setminus S$\/. 

On the other hand, notice that every $m\geq 2g$ is the $(m+1-g)$\/-th 
element of $S$\/, that is $\rho_{m+1-g}\,$, according to the notations 
of \cite{KirPel}. Finally, the semigroup $S$ is called {\em symmetric} 
when $r\in S$ if and only if $c-1-r\notin S$\/, for all $r\in\Z$. 
This is equivalent to say $c=2g$\/, that is to say $l_{g}=2g-1$.

\begin{defi}\label{FRdist}

For any semigroup $S\subseteq\N$ with $\sharp(\N\setminus S)<\infty$ 
and $0\in S$\/, the Feng-Rao distance of $S$ is defined by the function 
$$\delta_{FR}\;:\;S\longrightarrow\N$$
$$m\mapsto\delta_{FR}(m)\doteq\min\{\nu(r)\;|\;r\geq m,\;\;r\in S\}$$
where $\nu$ is the function 
$$\nu\;:\;S\longrightarrow\N$$
$$r\mapsto\nu(r)\doteq\sharp\{(a,b)\in S^{2}\;|\;a+b=r\}$$

\end{defi}

With the above notations, the following result summarizes some known facts 
about the functions $\nu$ and $\delta_{FR}$ for an arbitrary semigroup. 
One can check the details in \cite{HPNew} or \cite{KirPel}.

\begin{prop}\label{knownFR}

$\;$

\begin{description}

\item[(i)] $\delta_{FR}(0)=\nu(0)=1$, 
and $2\leq\nu(m),\delta_{FR}(m)\leq m+1$ if $m\in S\setminus\{0\}$\/. 

\item[(ii)] $\nu(m)=m+1-2g+D(m)$ for $m\geq c$ 
\footnote{Notice that one usually assumes that $m>2g-2$ in coding theory, 
and hence $m\geq c$\/. }, where 
$$D(m)\doteq\sharp\{(x,y)\;\;|\;\;x,y\;\;\mbox{are gaps of}\;\;S\;\;
\mbox{and}\;\;x+y=m\}$$

\item[(iii)] $\nu(m)=m+1-2g$ for all $m\in S$ with $m\geq 4g-1$. 

\item[(iv)] $\delta_{FR}(m)\geq m+1-2g\doteq d^{\ast}(m-1)$ for all $m\in S$\/, 
and equality holds if moreover $m\geq 4g-1$. 

\end{description}

\end{prop}

In particular, it follows that $\delta_{FR}(m)=\nu(m)=m+1-2g$ 
for all $m\in S$ such that $D(m)=0$. 
On the other hand, since from {\bf (iii)} one has that $\nu$ 
is an increasing function for $m\geq 4g-1$ at most, the knowledge 
of the Feng-Rao distance is finitely determined by the values of $\nu$\/.  
In fact, for any $m\in S$ one has 
$$\delta_{FR}(m)=\min\{\nu(m),\nu(m+1),\ldots,\nu(m')\}$$
where $m'$ is the least element in $S$ with $m'\geq m$ and 
such that $\nu(m')=m'+1-2g$\/. Thus, elements $m'\in S$ with 
$D(m')=0$ are interesting for calculations. Other interesting 
kind of elements in $S$ are those satisfying the following equality: 
$$\delta_{FR}(m)=\min\{r\in S\;|\;r\geq m+1-2g\}\hspace{1cm}[@]$$
In practice, such elements exist and have the advantage 
that the Feng-Rao distance for them is easy to compute. 

Our aim is to give some formulae which will allow us to compute the 
Feng-Rao distance for an arbitrary semigroup by means of an algorithm, 
provided a suitable set of generators is given. 
The main tool which will be used was introduced by Ap\'{e}ry in \cite{Ap} 
in order to study semigroups of curve singularities, and it is 
nothing but the Ap\'{e}ry systems of generators and their relations.

\begin{defi}\label{ApSet}

Let $S\subseteq\N$ be a semigroup with $\sharp(\N\setminus S)<\infty$ 
and $0\in S$\/; for $e\in S\setminus\{0\}$ define 
the Ap\'{e}ry set of $S$ related to $e$ by 
$$\{a_{0},a_{1},\ldots,a_{e-1}\}$$
where $a_{i}\doteq min\{m\in S\;|\;m\equiv i\;(mod\;\;e)\}$ 
for $0\leq i\leq e-1$. 

\end{defi}

Usually one takes $e$ as the {\em multiplicity} of $S$\/, that is, 
$e=e_{0}\doteq\min(S\setminus\{0\})$\/, but actually it is not necessary. 
On the other hand, notice that one could remove $a_{0}=0$ 
since it does not add any information about the semigroup. 
In the sequel, the index $i$ will be identified to the corresponding 
element in $\Z/(e)$\/. In fact, one has a disjoint union 
$$S=\bigcup_{i=0}^{e-1}(a_{i}+e\N)$$
and therefore the set $\{a_{1},\ldots,a_{e-1},e\}$ 
is a generator system for the semigroup $S$, called 
the Ap\'{e}ry (generator) system of $S$ related to $e$\/. 

Moreover, let $i,j\in\Z/(e)\equiv\Z_{e}$ 
and consider $i+j\in\Z_{e}\,$; then 
$$a_{i}+a_{j}=a_{i+j}+\alpha_{i,j}e$$
with $\alpha_{i,j}\geq 0$, by definition of the Ap\'{e}ry set. 
The numbers $\alpha_{i,j}$ are called {\em Ap\'{e}ry relations}\/. 

Under these conditions, every $m\in S$ can be written in a unique way as 
$m=a_{i}+le$, with $i\in\Z_{e}$ and $l\geq 0$. Thus, we can associate 
to $m$ two {\em Ap\'{e}ry coordinates}\/ $(i,l)\in\Z_{e}\times\N$. 

In order to compute $\nu(m)$\/, set $m\equiv(i,l)$\/, $a\equiv(i_{1},l_{1})$ 
and $b\equiv(i_{2},l_{2})$\/; since 
$$m=a+b=a_{i_{1}+i_{2}}+(l_{1}+l_{2}+\alpha_{i_{1},i_{2}})\,e$$
then $l_{1}+l_{2}=l-\alpha_{i_{1},i_{2}}\,$. 
Write $i_{1}=k$ and $i_{2}=i-k$\/; 
if $l<\alpha_{k,i-k}$ the equality $m=a+b$ is not possible, 
and so we are just interested in the case $\alpha_{k,i-k}\leq l$\/. 

Thus, for $0\leq i\leq e-1$ and $h\geq 0$ define 
$$B_{i}^{(h)}\doteq\sharp\{\alpha_{k,i-k}\leq h\;|\;k\in\Z_{e}\}$$
With these notations, the following result gives us a formula 
to compute $\nu(m)$\/.

\begin{prop}\label{formNm}

\hspace{.5cm}
$\nu(m)=B_{i}^{(0)}+B_{i}^{(1)}+\ldots+B_{i}^{(l)}$ 

\end{prop}

\dem

Suppose $\alpha_{k,i-k}=h\leq l$\/; 
then it has been considered at the right sum in the sets defining 
$B_{i}^{(h)},B_{i}^{(h+1)},\ldots,B_{i}^{(l)}$, that is $l-h+1$ times. 

On the other hand, the equality $l_{1}+l_{2}=l-\alpha_{k,i-k}$ holds 
for $l-h+1$ possible pairs $l_{1},l_{2}\,$, and so the theorem is proved. 

\vspace{.1cm}
\findemo

\newpage

Now, if we want to have a formula to compute the Feng-Rao distance, 
the main remark is that $\nu(m)$ is {\em increasing in $l$}\/, 
because of the previous formula. Then it suffices to calculate 
a {\em minimum in the coordinate} $i$\/, what gives only a finite number 
of possibilities. More precisely, one obtains the following result.

\begin{thm}\label{formFRd}

With the above notations, set $m=a_{i}+le$\/. 
For each $j\in\mbox\Z_{e}\,$, take $m_{j}=a_{j}+t_{j}e$\/, 
where $t_{j}$ is the minimum integer such that 
$t_{j}\geq max\left\{\displaystyle\frac{a_{i}-a_{j}}{e}+l,0\right\}$\/. 
Then one has 
$$\delta_{FR}(m)=min\{\nu(m_{j})\;|\;j\in\mbox\Z_{e}\}$$

\end{thm}

\dem

By using the above remark on the number $\nu(m)$\/, the formula 
follows from the fact that $m_{j}$ is the minimum element of $S$ 
with first Ap\'{e}ry coordinate equal to $j$ such that $m_{j}\geq m$\/. 

\vspace{.1cm}
\findemo

\vspace{.1cm}

As a conclusion, computing the Feng-Rao distance is easy if we have 
the Ap\'{e}ry set related to any non-zero element of the semigroup, and 
the method works in a quite general situation. 
Next, we will show how the above facts can be done more precise 
for the case of symmetric semigroups.

\subsection{Symmetric semigroups}

Now we will compute the value of $\delta_{FR}(m)$ for many 
elements $m\in S$\/, $S$ being a symmetric semigroup, and 
improve the computation in theorem \ref{formFRd} for them. 
The underlying idea is searching for the values $m\in S$ 
such that either $\delta_{FR}(m)=\nu(m)=d^{\ast}(m-1)=m+1-2g$ 
or the formula $[@]$ is satisfied; 
this formula will be called \lq\lq minimum formula" in the sequel. 
Our results will partially cover lacks in the results given in \cite{KirPel}. 

First of all, symmetry can be easily expressed in terms of Ap\'{e}ry sets. 
In fact, let $e$ be any non-zero element of the semigroup $S$\/, and 
consider the Ap\'{e}ry set $\{a_{0},a_{1},\ldots,a_{e-1}\}$ related to $e$\/. 
For this Ap\'{e}ry set, consider the index $N\in\Z_{e}$ such that 
$$a_{N}=\max\{a_{0},a_{1},\ldots,a_{e-1}\}$$
One can easily check that the last gap of $S$ is just $l_{g}=a_{N}-e$\/, 
i.e. $a_{N}=c-1+e$\/. 
Thus, it is easy to check that $S$ is symmetric if an only if 
$$a_{i}+a_{N-i}=a_{N}\hspace{1.cm}\forall i\in\Z_{e}$$
In this case, notice that $l_{g}=2g-1$, and thus $a_{N}=2g-1+e$\/. 
By using this fact and the formula of proposition \ref{formNm}, 
one obtains the following result, which provides us with a formula of 
$\delta_{FR}(m)$ for a certain range of values in $S$\/.

\begin{thm}\label{NmFRdSym}

Let $S$ be a symmetric semigroup. Then one has 
$$\delta_{FR}(m)=\nu(m)=m-l_{g}=m+1-2g=e$$
for all $m=2g-1+e$ with $e\in S\setminus\{0\}$\/. 

\end{thm}

\dem

Take an arbitrary non-zero $e\in S$\/, and let $N$ be the index in 
$\Z_{e}$ such that $a_{N}$ is the maximum Ap\'{e}ry element related to $e$\/. 
We will first prove the formula for all $m\in S$ such that 
$m\equiv N\;(mod\;\;e)$\/, i.e. for $m=a_{N}+le$\/, with $l\geq 0$. 

By proposition \ref{formNm}, one has 
$$\nu(m)=B_{N}^{(0)}+B_{N}^{(1)}+\ldots+B_{N}^{(l)}$$
where $B_{N}^{(h)}=\sharp\{\alpha_{k,N-k}\leq h\;|\;k\in\Z_{e}\}$\/. 
But because of the equality $a_{k}+a_{N-k}=a_{N}$ for all $k\in\Z_{e}\,$, 
one has $\alpha_{k,N-k}=0$ for all $k\in\Z_{e}\,$, and thus 
$B_{N}^{(h)}=e$ for all $h\geq 0$. Hence, $\nu(m)=(l+1)\,e=le+(a_{N}-l_{g})=
m-l_{g}\,$, and the formula is proved for $\nu(m)$\/. 
The formula for $\delta_{FR}(m)$ follows immediately 
from the definition of the Feng-Rao distance and 
proposition \ref{knownFR} {\bf (iv)}. 

In particular, for the fixed $e$ the statement holds for $l=0$, 
that is, for $m=a_{N}=l_{g}+e$\/. But $e$ is an arbitrary 
non-zero element in $S$\/, and thus the theorem is proved. 

\vspace{.1cm}
\findemo

Theorem \ref{NmFRdSym} shows in particular that the minimum formula 
$[@]$ holds for half of the elements in the interval $[c-1,2c-2]$\/. 
But, in fact, this conclusion of theorem \ref{NmFRdSym} is also true 
for many other elements in this interval. 
In order to show it, we forget for a moment of Ap\'{e}ry systems, and for 
$m\in I\doteq[c,2c-2]$ define two integers $n,q$ given respectively by the equalities 
$$m=n+c-1\hspace{.5cm},\hspace{.5cm}n+q=c-1$$
Then one has $m=2c-2-q$\/, and $n\in S$ iff $q\notin S$\/. 
Moreover, from proposition \ref{knownFR} {\bf (ii)} and the fact that 
$S$ is symmetric, it follows that 
$$\nu(m)=n+\nu(q)$$
where $\nu(q)\doteq 0$ for $q\notin S$\/. Thus, in particular, 
if $n\in S$ one recovers $\nu(m)=n$\/, what gives a different 
proof of theorem \ref{NmFRdSym}. 

Assume now that $n\notin S$ and let $n'$ be the least integer in $S$ 
such that $n'>n$\/. Consider the element $q'$ given by $n'+q'=c-1$ 
and let $\delta=\delta(q)=q-q'=n'-n$ be the distance of $q$ 
(resp. $n$\/) to a gap (resp. non-gap) of $S$\/. 
Notice that $\delta(q)\leq e_{0}-1$, $e_{0}$ being the multiplicity of $S$\/, 
since the interval $[n,n')$ consists of gaps and there are at most $e_{0}-1$ 
consecutive gaps in $S$\/.

\begin{thm}\label{improvedFR}

Let $S$ be a symmetric semigroup and $m\in I=[c,2c-2]$\/. 
With notations as above, one has 
$$\delta_{FR}(m)=\min\{n+\nu(q),n+1+\nu(q-1),\ldots,
n+\delta-3+\nu(q-\delta+3),n'\}=$$
$$=\min\{\nu(m),\nu(m+1),\ldots,\nu(m+\delta-3),n'\}$$
In particular, one has $\delta_{FR}(m)\leq n'$\/. 

\end{thm}

\dem

Since $n'\in S$ one has $\nu(c-1+n')=n'$\/; thus, 
from \ref{knownFR} {\bf (iv)} and the definition 
of the Feng-Rao distance one has 
$$\delta_{FR}(m)=\min\{\nu(m),\nu(m+1),\ldots,\nu(m+\delta)=n'\}$$
But $\nu(m+\delta-1)=n'-1+\nu(q'+1)\geq n'$\/, since 
$q'+1\in S$\/. Moreover, $\nu(m+\delta-2)=n'-2+\nu(q'+2)\geq n'$\/, 
since $q'+2\in S\setminus\{0\}$ and therefore $\nu(q'+2)\geq 2$ 
(notice that $q'$ can be negative but only in the case $n'=c$\/, 
and in this case $q'=-1$). This proves the theorem. 

\vspace{.1cm}
\findemo

\begin{nota}\label{CompImpSym}

According to the above theorem, the Feng-Rao distance can be computed 
by using the formula given in proposition \ref{formNm} and taking the 
minimum of only $\delta-1$ integers. In order to carry out this 
computation is better, in principle, to use the first equality 
in the formula of \ref{improvedFR}, since integers with the size 
of $q$ are smaller than those with the size of $m$\/. 
Next corollary complements the result of theorem \ref{NmFRdSym}. 

\end{nota}

\begin{coro}\label{corImpSym}

Let $S$ be a symmetric semigroup and $m\in I$\/. 
With notations as above, the minimum formula $[@]$ holds for $m$ 
if and only if one has 
$$\nu(\overline{q})\geq\delta(\overline{q})$$
for every $\overline{q}$ in the interval $[q'+3,q]$\/. 

In particular, this is satisfied for the cases $n\in S$\/, 
$\delta=\delta(q)=1$ or $2$, and when $\delta=3$ but $q$ 
is not an irreducible element of $S$ (i.e. a non-zero element 
of $S$ which is not a sum of two non-zero elements of $S$\/). 

\end{coro}

\dem

It follows from theorems \ref{NmFRdSym} and \ref{improvedFR}. 
Notice that $\nu(q)\geq 3$ if $\delta=3$ and $q$ is not irreducible. 

\vspace{.1cm}
\findemo

\begin{nota}\label{checkCor}

The condition in the above corollary is easy to check, in practice, 
for many values of $m$ as, for example, those with small values of $q$\/. 

\newpage

However, by using Ap\'{e}ry systems relative to any non-zero element $e$ 
of $S$ (for instance $e=e_{0}$\/) it is always possible to check it. 
In fact, $\nu(q)$ can be computed from proposition \ref{formNm} and, 
in order to compute $\delta=\delta(q)$\/, we can use the following result. 

\end{nota}

\begin{prop}\label{formDelta}

Assume that $q\in S$\/, with Ap\'{e}ry coordinates (related to $e$\/) 
given by $q=a_{i}+le$\/. Then $\delta$ is the least integer such that 
$q-\delta<a_{i-\delta}\,$. 

\end{prop}

\dem

On needs to find the value of $S$ such that $(q-\delta,q]\subset S$ 
and $q-\delta\notin S$\/. But if $q\equiv i\;(mod\;\;e)$ one has 
$q-j\equiv i-j\;(mod\;\;e)$ for all $j\in\Z_{e}\,$. Thus, 
the proposition follows from the definition of Ap\'{e}ry elements. 

\vspace{.1cm}
\findemo

Finally, let $q_{0}$ be the smallest integer in $S$ such that 
$\nu(q_{0})<\delta(q_{0})$\/, and set $q_{0}=c-1$ if such an 
integer does not exist. From corollary \ref{corImpSym} one has 
$q_{0}\geq e_{0}+2$ and $q_{0}=e_{0}+2$ if and only if 
$e_{0},e_{0}+1,e_{0}+2\in S$ and $e_{0}>2$. 
This happens, for instance, for the symmetric semigroup 
whose gaps are $1,2,\ldots,g-1,2g-1$, with $g\geq 3$. 
For a general symmetric semigroup, 
the value of $q_{0}$ is, in practice, much bigger. 

Now, from proposition \ref{knownFR} {\bf (iv)}, 
the minimum formula holds in the interval $[4g-1,\infty)$\/. 
In fact, from the above results one deduces the whole 
interval where this formula is satisfied.

\vspace{.3cm}

\begin{coro}\label{intervalSym}

For a symmetric semigroup $S$\/, the minimum formula $[@]$ 
is satisfied in the interval $(m_{0},\infty)$\/, where 
$m_{0}=4g-2-q_{0}\,$, but it does not hold for $m_{0}\,$. 

\end{coro}

\dem 

It follows from theorem \ref{improvedFR}, 
corollary \ref{corImpSym} and the definition of $q_{0}\,$. 

\vspace{.1cm}
\findemo

\begin{nota}

Notice that $m_{0}\geq 4g-e_{0}-3$, what gives an estimate for $m_{0}\,$. 

\end{nota}

\begin{ejplo}\label{exsSym}

$\;$

\begin{itemize}

\item For the semigroup $S=\langle 9,12,15,17,20,23,25,28 \rangle$ 
one has $c=32$ and $q_{0}=25$, thus $[@]$ is satisfied for $m\geq 38$. 

\item For the semigroup $S=\langle 6,8,10,17,19 \rangle$ 
one has $c=22$ and $q_{0}=19$, thus $[@]$ is satisfied for $m\geq 24$. 

\item For the semigroup $S=\langle 8,10,12,13 \rangle$ 
one has $c=28$ and $q_{0}=25$, thus $[@]$ is satisfied for $m\geq 31$. 

\newpage

\item For the semigroup $S=\langle 6,10,15 \rangle$ 
one has $c=30$ and $q_{0}=29$, thus $[@]$ is satisfied for $m\geq 30$. 

\end{itemize}

\end{ejplo}

\subsection{Semigroups at infinity}

In the situation of semigroups at infinity being computed as in 
section {\bf 3}, we want to calculate the Ap\'{e}ry set of the semigroup 
$S_{P}$ related to $e=\delta_{0}=deg\,\Upsilon$, 
with the assumptions of the Abhyankar-Moh theorem. In general, 
we can solve this problem for an arbitrary telescopic semigroup $S$\/. 

Recall that, for such a semigroup, the main fact is that 
any $m\in S$ can be written in a unique way in the form 
$$m=\sum_{k=0}^{h}\lambda_{k}\delta_{k}=\lambda_{0}\delta_{0}+
\sum_{k=1}^{h}\lambda_{k}\delta_{k}$$
with $\lambda_{0}\geq 0$ and $0\leq\lambda_{k}<n_{k}=d_{k}/d_{k+1}$ 
for $1\leq k\leq h$\/. 
But this means exactly that the possible {\em Ap\'{e}ry elements}\/ 
related to $\delta_{0}$ are only those with $\lambda_{0}=0$\/. 
On the other hand, the number of all the possible elements of this form is 
$$n_{1}n_{2}\ldots n_{h}=\delta_{0}$$
since $d_{1}=\delta_{0}$ and $d_{h+1}=1$. Thus, all of them are different 
modulo $\delta_{0}$ and minimum with this property, and hence they are 
just the elements of the Ap\'{e}ry set of $S$ related to $\delta_{0}\,$. 

As a consequence, for a given $m\in S$ written in the form 
$$m=\lambda_{0}\delta_{0}+\sum_{k=1}^{h}\lambda_{k}\delta_{k}$$
we can compute its {\em Ap\'{e}ry coordinates}\/ $l=\lambda_{0}\,$
and $i\in\Z/(\delta_{0})$ such that 
$i\equiv\sum_{k=1}^{h}\lambda_{k}\delta_{k}\;(mod\;\;\delta_{0})$\/. 
Thus, we can easily compute $a_{i}$ and 
$\alpha_{i,j}\,$, and so $\nu(m)$ and $\delta_{FR}(m)$ for a telescopic 
semigroup and, in particular, for a semigroup at infinity $S_{P}$ when 
computed by the algorithm of approximate roots, under the assumptions of 
the Abhyankar-Moh theorem.

In fact, the values of $\delta_{FR}(m)$ are known for many values of 
$m\in S$\/, if $S$ is telescopic. More precisely, apart from the 
results given by proposition \ref{knownFR}, if we assume moreover that 
$\delta_{h}=\max\{\delta_{0},\delta_{1},\ldots,\delta_{h}\}$\/, 
then the minimum formula $[@]$ is true for every $m\in[c,2c-2]$ 
in the interval \footnote{This result is improved by corollary \ref{intervalSym}. 
For example, the telescopic semigroup $S=\langle 8,10,12,13 \rangle$ 
the formula $[@]$ is satisfied from $m=31$ instead of $m=42$, 
and for the telescopic semigroup $S=\langle 6,10,15 \rangle$ 
the formula $[@]$ is satisfied from $m=30$ instead of $m=44$. } 
$[(4g-1)-(d_{h}-1)\delta_{h},4g-1)$\/, i.e. 
one has $q_{0}\geq(d_{h}-1)\delta_{h}\,$. Also one has $\delta_{FR}(m)=j+1$ 
if $(j-1)\delta_{h}<m\leq j\,\delta_{h}\leq(d_{h}-1)\delta_{h}\,$ 
(see \cite{KirPel}). Hence, the Feng-Rao distance is unknown in general 
for the interval $[(d_{h}-1)\delta_{h}+1,(4g-2)-(d_{h}-1)\delta_{h}]$\/, 
but nevertheless it can be computed by the formulae of proposition 
\ref{formNm} and theorem \ref{formFRd}. 
Notice that, for example, one has 
$$\nu\left(\Sum_{k=1}^{h}\lambda_{k}\delta_{k}\right)=
(\lambda_{1}-1)(\lambda_{2}-1)\ldots(\lambda_{h}-1)$$
for the Ap\'{e}ry elements.

Now, we want to compute the Ap\'{e}ry set related to $e=deg\,\Upsilon$ for 
the Weierstrass semigroup $\Gamma_{P}=S_{P}+b_{1}\N+\ldots+b_{s}\N$, 
where $S_{P}$ has been computed by means of the Abhyankar-Moh theorem, 
$s=dim_{\Ff}(A/B)$ and $b_{i}$ has been computed as in proposition \ref{TriangLem}. 
In order to do it, it suffices to solve the following question: 

\begin{quote}

For a given numerical semigroup $S$ with Ap\'{e}ry set 
$\{a_{0},a_{1},\ldots,a_{e-1}\}$ related to $e\in S\setminus\{0\}$\/, 
computing the Ap\'{e}ry set 
$\{\overline{a_{0}},\overline{a_{1}},\ldots,\overline{a_{e-1}}\}$ 
for the semigroup $\overline{S}=S+b\N$ related to $e\in\overline{S}\setminus\{0\}$\/. 

\end{quote}

\noindent
In fact, after having solved the above question, 
by repeating the procedure a finite number of steps one gets 
the Weierstrass semigroup $\Gamma_{P}$ with its Ap\'{e}ry basis 
(see remark \ref{fastTA}).

First, for a given $b$ one wants to know wether $b$ is already in $S$ or not, 
because if it is so, then $S=S+b\N$ and nothing changes. 
This can be easily tested with aid of the given Ap\'{e}ry set of $S$ as follows: 

\begin{description}

\item[(i)] Calculate $i\in\Z_{e}$ such that $i\equiv b\;(mod\;\;e)$\/. 

\item[(ii)] Then, $b\in S$ if and only if $l=\displaystyle\frac{b-a_{i}}{e}\geq 0$. 

\end{description}

On the other hand, 
the candidates to be the elements of the new Ap\'{e}ry set 
are obviously the numbers $m_{j,\lambda}=a_{j}+\lambda\,b$ 
with $0\leq j\leq e-1$ and $0\leq\lambda\leq e-1$, i.e. 
the number of possibilities is small. 
Thus, we proceed as follows: 

\begin{description}

\item[1. ] Initialize $\overline{a_{i}}=a_{i}$ for $1\leq i\leq e-1$ 
(obviously $\overline{a_{0}}=a_{0}=0$). 

\item[2. ] Take one of the elements $m_{j,\lambda}$ and compute 
its remainder $i$ modulo $e$\/. 

\item[3. ] Compare $m_{j,\lambda}$ with the value of $\overline{a_{i}}$\/; 
if $\overline{a_{i}}$ is greater we have to change its value, set 
$\overline{a_{i}}=m_{j,\lambda}$ and so on with the next $m_{j,\lambda}\,$. 

\end{description}

\noindent
At the end of this procedure we get the Ap\'{e}ry set for $\overline{S}$\/, 
and the new Ap\'{e}ry relations can be easily computed from them.

\begin{nota}\label{sGsAp}

The results of this section suggest to present the elements $m$ 
in a semigroup $S$ as couples $(i,l)\in\Z_{e}\times\N$ if $m=a_{i}+l\,e$. 
Thus, $S$ is given by nothing but the data 
$\{e \; ; \; a_{i} \;\,{\rm for}\;\, 0\leq i\leq e-1\}$ 
with the restrictions $\alpha_{i,j}\geq 0$. In particular, 
the computations made in section {\bf 3.3} to find a function 
for every element in the Weierstrass semigroup are easier 
if we use the Ap\'{e}ry description of the semigroup and functions 
for each element of the Ap\'{e}ry set. 

More precisely, one just saves the data 
$\{a_{0}=0,a_{1},\ldots,a_{e-1},e\}$ 
together with the associated functions 
$\{h_{0}=1,h_{1},\ldots,h_{e-1},h_{e}\}$ satisfying 
$-\upsilon_{P}(h_{i})=a_{i}$ for $0\leq i \leq e-1$, 
and $-\upsilon_{P}(h_{e})=e$\/. 
Thus, for any $m\in S$ we compute its coordinates $(i,l)$ 
and then $f_{m}=h_{i}\cdot h_{e}^{l}$ is the 
function associated to $m$\/. We save in general many generators, 
but the description of an element $m\in S$ in terms of them is 
quite fast (just doing a sum and two divisions by $e$\/). 
However, if we save the sequence of non-gaps $\{\rho_{i}\}$ with 
the associated functions $\{f_{i}\}$ we may use much more space, 
since the ratio $g/e$ can be as large as we want. 

On the other hand, the new elements added in the steps of the algorithm 
\ref{TriangAlg} modify the semigroup changing those data. Thus, such changes 
can be included in this algorithm, and so the semigroup $\Gamma_{P}\,$, 
the functions and the Feng-Rao distance can be computed simultaneously. 
In particular, this gives an effective answer to the idea which was proposed 
in remark \ref{fastTA}. 

Finally, for the computation of the Feng-Rao distance 
we should take into account that in case of $\Gamma_{P}$ 
being symmetric\footnote{With the notations of section {\bf 3.1}, 
this is equivalent to ${\cal C}$ being an affine complete intersection curve, 
as you can see in \cite{Sat}. }, then the calculations can be simplified 
by means of theorem \ref{improvedFR} or its corollaries. 

\end{nota}

\end{document}